\documentclass{article}%
\usepackage{amsfonts}
\usepackage{amsmath}
\usepackage{amssymb}
\usepackage{amsxtra}
\usepackage{graphicx}
\usepackage{geometry}
\usepackage{color}
\usepackage{colortbl}
\usepackage{caption}
\usepackage[draft]{varioref}%
\providecommand{\U}[1]{\protect\rule{.1in}{.1in}}
\newtheorem{theorem}{Theorem}

\newtheorem{definition}[theorem]{Definition}

\newtheorem{lemma}[theorem]{Lemma}

\newtheorem{proposition}[theorem]{Proposition}
\newtheorem{remark}[theorem]{Remark}

\numberwithin{equation}{section}

\begin{document}

\title{Bingham flow in porous media with obstacles of different size}
\author{R.Bunoiu\thanks{Institut Elie Cartan de Lorraine, CNRS, UMR 7502, Universit\'{e} de Lorraine, Metz, F-57045, France. email: renata.bunoiu@univ-lorraine.fr.},  G.Cardone\thanks{Universit\`{a} del Sannio, Department of Engineering, Corso Garibaldi, 107, 82100 Benevento, Italy; member of GNAMPA (INDAM). email: giuseppe.cardone@unisannio.it}}
\maketitle

\begin{abstract}
By using the unfolding operators for periodic homogenization, we give a
general compactness result for a class of functions defined on bounded domains
presenting perforations of two different size. Then we apply this result to
the homogenization of the flow of a Bingham fluid in a porous medium with
solid obstacles of different size. Next we give the interpretation of the
limit problem in term of a non linear Darcy law.

\medskip

Keywords: homogenization, unfolding operators, Bingham fluid, porous media

\medskip

MSC: 35B27, 76M50, 76S05

\end{abstract}

\section{Introduction\label{section1}}

In this paper we study the homogenization problem for a Bingham flow in a
porous medium with solid obstacles of different size. The aim of our paper is
twofold: we first define the unfolding operators for periodic homogenization
in a domain which presents periodically distributed perforations of two
different size and we give corresponding compactness results. Then we
illustrate these results with an application to the homogenization of a
Bingham flow in a porous medium with solid obstacles of different size.

In order to define the appropriate unfolding operators and to get the
compactness results, we follow the ideas introduced by D. Cioranescu, A.
Damlamian and G. Griso in \cite{4} and \cite{8} for the case of functions with one scale of
periodicity and developed later by A. Damlamian, N. Meunier, J. Van
Schaftingen in \cite{DaMu-VS} and \cite{Mu-VS} for the case of functions with
more than one periodicity scales. Nevertheless, our result is different from
the ones presented in the previous cited papers, due to the presence of the
perforations at the two different scales. The case, different from  the one presented here, corresponding
to the unfolding operators
for  a doubly periodic domain presenting perforations at the very small  scale only, was recently  addressed by Bunoiu, and Donato in \cite{BP}.

More precisely, our domain contains small perforations of size $\varepsilon$
periodically distributed with period $\varepsilon$ and very small perforations
of size $\varepsilon\delta(\varepsilon)$ periodically distributed with
periodicity $\varepsilon\delta(\varepsilon)$. Here $\varepsilon$ and
$\delta(\varepsilon)$ are real positive parameters smaller than one with
$\delta(\varepsilon)$ tending to zero when $\varepsilon$ tends to zero. Such a
geometry modelizes, for example, a porous medium in which the perforations
correspond to solid impervious obstacles.

In the fluid part of this porous medium we consider the stationary flow of the
Bingham fluid, under the action of external forces. The Bingham fluid is an
incompressible fluid which has a non linear constitutive law; so it is a
non-Newtonian fluid. This fluid moves like a rigid body when a certain
function of the stress tensor is below a given threshold. Beyond this
threshold the fluid flows, obeying a non linear constitutive law. As example of such fluids we can mention some paints, the mud which can be
used for the oil extraction and the volcanic lava. Bingham flow in other contexts is studied by  Bunoiu, and Kesavan in\cite{2}. For a presentation of the different types of non-Newtonian fluids we refer the reader to \cite{CioGiRa}. 

The mathematical model of the Bingham flow in a bounded domain was introduced
in \cite{6} by G. Duvaut and J. L. Lions. The existence of the velocity and of
the pressure for this model was proved in the case of a bi-dimensional and of
a three-dimensional domain.

The homogenization problem in a classical porous medium, with obstacles of
size $\varepsilon$ and $\varepsilon$-periodically distributed, was first
studied in \cite{7} by J. L. Lions and E. Sanchez-Palencia. The authors did
the asymptotic study of the problem by using a multiscale method, involving a
\textquotedblleft macroscopic\textquotedblright\ variable $x$ and a
\textquotedblleft microscopic\textquotedblright\ variable $y=\dfrac
{x}{\varepsilon}$, associated to the relative dimension of the pores. The
study is based on a multiscale \textquotedblleft ansatz\textquotedblright,
which allows to obtain to the limit a non linear Darcy law. There is no
convergence result proved.

The rigorous justification for the convergence of the homogenization process
of the results presented in \cite{7} is given by A. Bourgeat and A. Mikelic in
\cite{1}. In order to do it, the authors used monotonicity methods coupled
with the two-scale convergence method introduced by G. Nguetseng in \cite{10}
and further developed by G. Allaire in a series of papers, as for example
\cite{9}. The limit problem announced in \cite{7} was obtained, by letting the
small parameter $\varepsilon$ tend to zero in the initial problem. The
unfolding method for periodic homogenization, introduced by D. Cioranescu, A.
Damlamian and G. Griso in \cite{4} was used by R. Bunoiu, G. Cardone, C.
Perugia in \cite{BuCaPe} in order to obtain the limit problem. This method
presents the advange of transforming in an easy manner the initial problem,
stated in a domain dependent on $\varepsilon$, in a problem stated in a domain
independent of $\varepsilon$. The passage to the limit when $\varepsilon$
tends to zero is then simple thanks to the compactness results, and this for
the non linear terms too.

Our paper is organized as follows. In section \ref{section2} we define the
double perforated domain, the unfolding operators adapted to it and then we
give a compactness result. In section \ref{section3} we describe the problem
of the Bingham flow and we give the \textit{a priori} estimates for the
velocity and the pressure of the flow. Following the ideas in R. Bunoiu, J.
Saint Jean Paulin \cite{BuSJP}, we construct the extension of the pressure to
the whole domain, namely in the perforations too. In section \ref{section4} we
state the main result of the paper, which consists in getting the limit
problem. It is obtained in two steps: we first apply the unfolding operators
for periodic homogenization defined in section \ref{section2} to the
variational formulation of the problem which describes the Bingham flow in our
porous medium. Then we pass to the limit when the small parameter
$\varepsilon$ tends to zero. In section \ref{section5} we give the
interpretation of the limit problem in term of a non linear Darcy law and we
compare it with the classical linear Darcy law.

\section{Unfolding operators and compactness results\label{section2}}

Let $\Omega$ be a bounded open domain in $\mathbb{R}^{n}$ with Lipschitz
continuous boundary $\partial\Omega,$ $n=2$ or $n=3$.

We consider two fixed reference cells $Y=\left]  0,y_{1}^{0}\right[
\times...\times\left]  0,y_{n}^{0}\right[  $ and $Z=\left]  0,z_{1}%
^{0}\right[  \times...\times\left]  0,z_{n}^{0}\right[  $ and two closed
subsets $Y_{s}$ and $Z_{s}$ with non-empty interior and Lipschitz continuous
boundaries, contained in $Y$ and $Z$ respectively. We define:%
\[
Y^{\ast}=Y\setminus Y_{s},\ \ Z^{\ast}=Z\setminus Z_{s},
\]
and we give a simple example in Figure 1.%

\begin{figure}
\begin{center}
\includegraphics[scale=0.65]{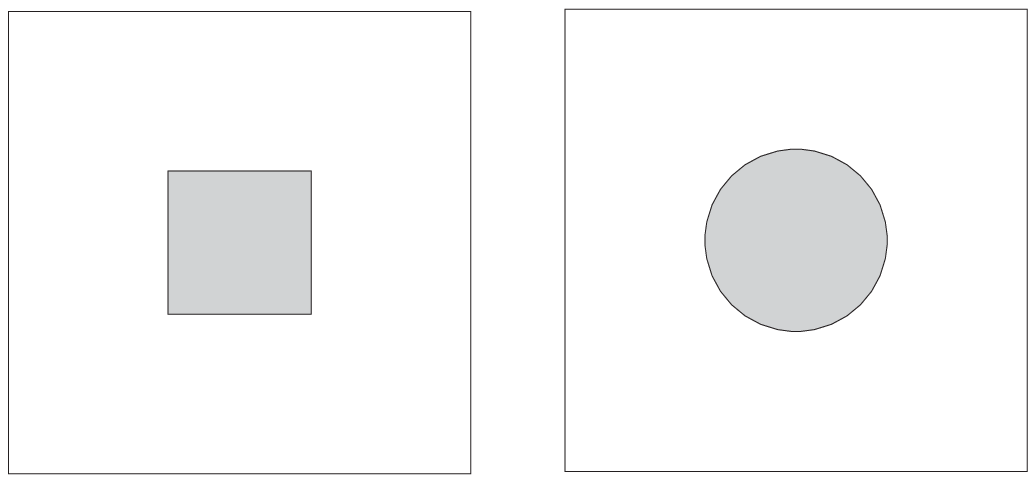}
\end{center}
\caption{Domains $Y^{\ast}$ and $Z^{\ast}$}
\label{f1}
\end{figure}

Let $\varepsilon$ be a positive parameter, smaller than one. For every
$\varepsilon>0$, let $0<\delta(\varepsilon)<\varepsilon$ be such that
\[
\lim\limits_{\varepsilon\rightarrow0}\delta(\varepsilon)=0.
\]
We suppose that there exists an $\varepsilon$ such that the domain
$\overline{Y}$ is exactly covered by a finite number of cells $\delta
(\varepsilon)\overline{Z}$. Moreover, we suppose that $Y_{s}$ is exactly
covered by a finite number of cells $\delta(\varepsilon)\overline{Z}$. This
last hypothesis implies some restrictions for the geometry of $Y_{s}$. We
deduce that there is no intersection between the domains $Y_{s}$ and
$\delta(\varepsilon)Z_{s}$ in the cell $Y$, as one can see on an example in
Figure 2. If we consider all the small parameters $\varepsilon/2^{N}$ (with
$N$ natural number), the above assumptions are still true. We denote $Y_{f}$
the complement in $Y$ of the set ${Y_{s}}\bigcup\limits_{l\in{\mathbb{Z}}^{n}%
}({\delta(\varepsilon)(lz^{0}+{Z}_{s})}\cap Y^{\ast})$.%

\begin{figure}
\begin{center}
\includegraphics[scale=0.65]{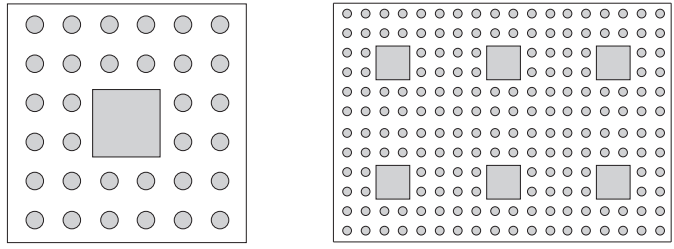}
\end{center}
\caption{Domains $Y_{f}$ and $\Omega_{\varepsilon\delta}$}
\label{f2}
\end{figure}

We multiply the perforated cell $Y$ (Figure 2) by $\varepsilon$ and we repeat
it in the domain $\Omega$. For simplicity and without loosing any generality,
one could even assume that $\overline{\Omega}$ is exactly covered by a finite
number of cells $\varepsilon\overline{Y}$. The domain $\Omega_{\varepsilon
\delta}$ is the one obtained by taking out of $\Omega$ the translated of the
domains $\varepsilon Y_{s}$ and $\varepsilon\delta(\varepsilon)Z_{s}$. Let us
notice that there is no intersection between the solid obstacles $\varepsilon
Y_{s}$ and $\varepsilon\delta(\varepsilon)Z_{s}$ in $\Omega_{\varepsilon
\delta}$, because there is no intersection between the solid obstacles $Y_{s}$
and $\delta(\varepsilon)Z_{s}$ in the cell $Y$. The domain $\Omega
_{\varepsilon\delta}$ is connected, but the union of solid obstacles is not
connected (see an example in Figure 2).

Let $\chi_{Y^{\ast}}$ and $\chi_{Z^{\ast}}$ be the characteristic functions of
the domains $Y^{\ast}$and $Z^{\ast}$, defined by:%
\[
\chi_{Y^{\ast}}\left(  y\right)  =\left\{
\begin{array}
[c]{l}%
1,\ \text{in }Y^{\ast},\\
0,\ \text{in }Y\setminus Y^{\ast}%
\end{array}
\right.  ,\ \ \ \chi_{Z^{\ast}}\left(  y\right)  =\left\{
\begin{array}
[c]{l}%
1,\ \text{in }Z^{\ast},\\
0,\ \text{in }Z\setminus Z^{\ast}.
\end{array}
\right.
\]
We extend the characteristic functions $\chi_{Y^{\ast}}$ (respectively
$\chi_{Z^{\ast}}$) by periodicity, with period $y_{i}^{0}$ in $y_{i}$ and with
period $z_{i}^{0}$ in $z_{i}$, for $i=1,..,n$. The domain $\Omega
_{\varepsilon\delta}$, defined as above is described by:%
\[
\Omega_{\varepsilon\delta}=\left\{  x\in\Omega:\chi_{Y^{\ast}}\left(  \frac
{x}{\varepsilon}\right)  \chi_{Z^{\ast}}\left(  \frac{x}{\varepsilon
\delta(\varepsilon)}\right)  =1\right\}  .
\]
The domain $\Omega_{\varepsilon\delta}$ presents a structure with a double
periodicity: there are small perforations of size $\varepsilon$ and very small
perforations of size $\varepsilon\delta(\varepsilon)$. The boundary
$\partial\Omega_{\varepsilon\delta}$ that is composed by two parts: the
boundary of $\Omega$, denoted $\partial\Omega$ and the union of the boundaries
of all the obstacles, denoted $\Gamma_{\varepsilon}.$

\bigskip

We follow the general idea of the unfolding method, namely we transform
oscillating functions defined on the domain $\Omega$ into functions defined on
the domain $\Omega\times Y\times Z$. In order to do this, we proceed in two
steps: first we use the general theory of the unfolding homogenization in
order to make the transformation from the domain $\Omega$ to the domain
$\Omega\times Y$. In order to do this, we use the unfolding operator
introduced in \cite{4} for the scale $\varepsilon$. Next we define a second
unfolding operator, for the scale $\delta(\varepsilon)$, which allows us to
transform oscillating functions defined on $\Omega\times Y$ into functions
defined on $\Omega\times Y\times Z$. In order to do this, we will follow in
addition the ideas of A. Damlamain N. Meunier, J. Van Schaftingen in
\cite{DaMu-VS} and \cite{Mu-VS}.

For the first step, the idea is to transform oscillating functions defined on
the domain $\Omega$ into functions defined on the domain $\Omega\times Y$, in
order to isolate the oscillations in the second variable. This transformation,
together with \textit{a priori} estimates, allows us to use compactness
results and then to get the limits of our oscillating sequences. We start by
recalling the results as far as the unfolding operator for the scale
$\varepsilon$ is concerned.

\bigskip

We know that every real number $a$ can be written as the sum between his
integer part $[a]$ and his fractionary part $\left\{  a\right\}  $ which
belongs to the interval $\left[  0,1\right)  $.

For $x=(x_{1},...,x_{n})\in\mathbb{R}^{n}$, we apply a similar decomposition
to every real number $\dfrac{x_{i}}{\varepsilon}$ for $i=1,...,n$ and we get
\[
x=\varepsilon\Bigl[\frac{x}{\varepsilon}\Bigr]_{Y}+\varepsilon\Bigl\{\frac
{x}{\varepsilon}\Bigr\}_{Y},
\]
where $\Bigl[\dfrac{x}{\varepsilon}\Bigr]_{Y}\in\mathbb{Z}^{n}$ and
$\Bigl\{\dfrac{x}{\varepsilon}\Bigr\}_{Y}\,\in Y.$

We define now
\[
\widehat{K}_{\varepsilon}=\{ k\in\mathbb{Z}^{N} \, \vert\, \varepsilon(k+
Y)\subset\Omega\}, \quad\widehat\Omega_{\varepsilon}= \mathrm{int}%
\bigcup\limits_{k\in\widehat{K}_{\varepsilon}}\left(  \varepsilon(k+
\overline{Y})\right)  , \quad\Lambda_{\varepsilon}=\Omega\setminus
\widehat\Omega_{\varepsilon},
\]

and we notice that the set $\widehat\Omega_{\varepsilon}$ is the interior of
the largest union of $\varepsilon(k+ \overline{Y}) $ cells included in
$\Omega$.

\begin{definition}
For any Lebesgue measurable function $\varphi$ on $\Omega$, we define the
periodic unfolding operator by the formula
\[
T_{\varepsilon}(\varphi)(x,y)=\left\{
\begin{array}
[c]{ll}%
\varphi\left(  \varepsilon\left[  \dfrac{x}{\varepsilon}\right]  _{Y}+
\varepsilon y \right)  & \text{for a.e. } (x,y)\in\widehat\Omega_{\varepsilon
}\times Y,\\[2ex]%
0 & \text{for a.e. } (x,y)\in\Lambda_{\varepsilon}\times Y.
\end{array}
\right.
\]
\label{T-epsilon}
\end{definition}

According to \cite{4}, this operator has the following properties:

\begin{enumerate}
\item[p$_{1})$] $T_{\varepsilon}$ is linear and continuous from $L^{2}%
(\Omega)$ to $L^{2}(\Omega\times Y)$;

\item[p$_{2})$] $T_{\varepsilon}(\varphi\,\phi)=T_{\varepsilon}(\varphi
)T_{\varepsilon}(\phi),$ $\forall\varphi,\phi\in L^{2}(\Omega);$

\item[p$_{3})$] If $\varphi\in L^{2}(Y)$ is a $Y$-periodic function and
$\varphi^{\varepsilon}(x)=\varphi\Bigl(\displaystyle\frac{x}{\varepsilon
}\Bigr),\,x\,\in\,\mathbb{R}^{N}$ then%
\[
T_{\varepsilon}(\varphi_{|_{\Omega}}^{\varepsilon})\,\rightarrow\varphi\text{
strongly in }L^{2}(\Omega\times Y);
\]

\item[p$_{4})$] If $\varphi_{\varepsilon}\in L^{2}(\Omega)$ and $\varphi
_{\varepsilon}\rightarrow\varphi$ strongly in $L^{2}(\Omega),$ then%
\[
T_{\varepsilon}(\varphi_{\varepsilon})\rightarrow\varphi\text{ strongly in
}L^{2}(\Omega\times Y).
\]

\end{enumerate}

Moreover, the following results hold (see Proposition 2.9 (iii) in \cite{4}):

\begin{proposition}
\label{prop1} Let $\{\varphi_{\varepsilon}\}_{\varepsilon}$ be a bounded
sequence in $L^{2}(\Omega)$ such that
\[
T_{\varepsilon}(\varphi_{\varepsilon})\to{\varphi}\text{ weakly in }%
L^{2}(\Omega\times Y).
\]
Then%
\[
{\varphi}_{\varepsilon}\to\mathcal{M}_{Y}({\varphi})\text{ weakly in }%
L^{2}(\Omega),
\]
\text{where the mean value operator} $\mathcal{M}_{Y}$ \text{is defined by}
\[
\mathcal{M}_{Y}({\varphi})= \frac{1}{\left\vert Y\right\vert }\int_{Y}%
{\varphi}\left(  x,y\right)  dy \text{ a.e. for } x \in\Omega.
\]

\end{proposition}

\bigskip

\noindent Let us moreover observe that for a function $\varphi\,\in
H^{1}(\Omega)$, one has%
\[
\nabla_{y}(T_{\varepsilon}(\varphi))=\varepsilon T_{\varepsilon}(\nabla
\varphi)\quad\text{ a.e. }(x,y)\,\in\Omega\times Y.
\]

We define now the second unfolding operator, at the scale $\delta
(\varepsilon)$ (denoted in the sequel by $\delta$).

\begin{definition}
\label{T-delta} Let $v \in L^{p}(\Omega\times Y)$. Then the unfolding operator
$T_{\delta}$ is defined by%

\[
T_{\delta}(v)(x,y,z)=v \left(  x, \delta\Bigl[\frac{y}{\delta}\Bigr]_{Z} +
\delta z \right)  \quad\text{ for } (x,y,z) \in\Omega\times Y \times Z,
\]
where $x$ plays the role of a parameter.
\end{definition}

Let now $\varPhi$ be a function belonging to the space $H^{1}_{0}(\Omega)$.
Then, accordind to Definitions \ref{T-epsilon} and \ref{T-delta} we have%

\[
T_{\delta}(T_{\varepsilon}(\varPhi))(x,y,z)= \varPhi\Bigl(\varepsilon
\Bigl[\frac{x}{\varepsilon} \Bigr]_{Y} +\varepsilon\delta\Bigl[\frac{y}%
{\delta} \Bigr]_{Z} + \varepsilon\delta z \Bigr) .
\]

Moreover, the following equality holds true:%

\[
\varepsilon\delta\displaystyle T_{\delta}(T_{\varepsilon}( \nabla_{x}
\varPhi))= \displaystyle T_{\delta}(T_{\varepsilon}(\varepsilon\delta
\nabla_{x} \varPhi))= T_{\delta}(\delta T_{\varepsilon}(\varepsilon\nabla_{x}
\varPhi_{\varepsilon\delta} ))= T_{\delta}(\delta(\nabla_{y} T_{\varepsilon
}(\varPhi))= \nabla_{z} T_{\delta}(T_{\varepsilon}(\varPhi))
\]

and we have the convergence results:

\begin{proposition}
\label{prop2} Let $\left\{  \varPhi_{\varepsilon\delta}\right\}
_{\varepsilon\delta}$ be a sequence in $H_{0}^{1}(\Omega)$ bounded in
$L^{2}(\Omega)$. Let us assume that%
\[
\varepsilon\delta\left\Vert \nabla\varPhi_{\varepsilon\delta}\right\Vert
_{(L^{2}\left(  \Omega\right)  )^{n}}\leq C.
\]
Then there exists $\widehat{\varPhi}$ in $L^{2}(\Omega\times Y;H^{1}%
_{\#}\left(  Z\right)  )$ such that, up to a subsequence still denoted by
$\varepsilon\delta$ we have
\begin{align*}
T_{\delta}\left(  T_{\varepsilon}(\varPhi_{\varepsilon\delta})\right)   &
\rightarrow\widehat{\varPhi}\text{ weakly in }L^{2}(\Omega\times Y\times Z),\\
T_{\varepsilon}(\varPhi_{\varepsilon\delta})  &  \rightarrow\frac{1}{|Z|}%
\int_{Z}\widehat{\varPhi}dz\text{ weakly in }L^{2}(\Omega\times Y),\\
\varPhi_{\varepsilon\delta}  &  \rightarrow\frac{1}{|Y||Z|}\int_{Y\times Z}\widehat{\varPhi}dydz\text{ weakly in }L^{2}(\Omega),\\
\varepsilon\delta T_{\delta}\left(  T_{\varepsilon}(\nabla_{x}%
\varPhi_{\varepsilon\delta})\right)   &  \rightarrow\nabla_{z}\widehat
{\varPhi}\text{ weakly in }(L^{2}(\Omega\times Y\times Z))^{n},
\end{align*}
where $\displaystyle H^{1}_{\#}\left(  Z \right)  =\left\{  \phi\in H^{1}(Z),
\phi\hbox { is Z-periodic} \right\}  .$
\end{proposition}

\textbf{Proof.} The sequence $\left\{  \varPhi_{\varepsilon\delta}\right\}
_{\varepsilon\delta}$ being bounded in $H_{0}^{1}(\Omega)$, then $\left\{
T_{\delta}\left(  T_{\varepsilon}(\varPhi_{\varepsilon\delta})\right)
\right\}  _{\varepsilon\delta}$ is bounded in the space $L^{2}(\Omega\times
Y\times Z)$. Clasical compactness results imply the existence of a function
$\widehat{\varPhi}$ in $L^{2}(\Omega\times Y\times Z)$ such that the first
weak convergence holds true. By using Proposition \ref{prop1} and its
analoguous at the scale $\delta$, we obtain the second and the third weak convergences.

The last weak converge is a consequence of the equality
\[
\varepsilon\delta\displaystyle T_{\delta}(T_{\varepsilon}(\nabla
_{x}\varPhi_{\varepsilon\delta}))=\nabla_{z}T_{\delta}(T_{\varepsilon
}(\varPhi_{\varepsilon\delta}))
\]
and of the first weak convergence, for the sequence $T_{\delta}(T_{\varepsilon
}(\nabla_{x}\varPhi_{\varepsilon\delta})).$ The fact that the limit
$\widehat{\varPhi}$ actually belongs to the space $L^{2}(\Omega\times
Y;H^{1}_{\#}\left(  Z\right)  )$ is due to the application of a result from
\cite{4} to the unfolding operator at the scale $\delta$. $\blacksquare$

\begin{remark}
\label{rk} If $\left\{  \varPhi_{\varepsilon\delta}\right\}  _{\varepsilon
\delta}$ is a sequence in $H_{0}^{1}(\Omega_{\varepsilon\delta})$ those
extension by zero to the whole of $\Omega$ satisfy the hypothesis of
Proposition \ref{prop2}, then all the results still hold true, with $Y$ and
$Z$ replaced by $Y^{\ast}$ and $Z^{\ast}$ respectively in the description of
all the function spaces and in the integrals. Indeed, the sequence $\left\{
\varPhi_{\varepsilon\delta}\right\}  _{\varepsilon\delta} $ vanishes on
$\Omega\times Y_{s} \times Z_{s}$ and this property is preserved by passing to
the limit.
\end{remark}

To end this section, we recall one of the key points of the use of the
unfolding method for periodic homogenization: the fact that the integrals over
the domain $\Omega$ can be replaced by integrals over the domain $\Omega\times
Y\times Z$, by using the relation below
\[
\int_{\Omega}\varPhi dx\sim\dfrac{1}{\left\vert Y\right\vert \left\vert
Z\right\vert }\int_{\Omega\times Y\times Z}T_{\delta}(T_{\varepsilon}\left(
\varPhi\right)  )dxdydz,\ \ \forall\varphi\in L^{1}\left(  \Omega\right)  ,
\]
which is true for $\varepsilon$ and $\delta$ sufficiently small.

\section{Statement of the problem and preliminary results\label{section3}}

Our aim now is to apply the results from Section \ref{section2} to the
homogenization of a problem stated in a domain $\Omega_{\varepsilon\delta}$
defined as before. The problem we address is the flow of a Bingham fluid in a
porous medium with obstacles of different size. Indeed, such porous media can
be modelized by the domain $\Omega_{\varepsilon\delta}$, corresponding to the
part where the fluid flows. The perforations correspond to solid impervious
obstacles. If $u_{\varepsilon\delta}$ and $p_{\varepsilon\delta}$ are the
velocity and pressure respectively for a Bingham fluid, then its stress tensor
is defined by%

\begin{equation}
\sigma_{ij}=-p_{\varepsilon\delta}\delta_{ij}+g \varepsilon\delta\dfrac
{D_{ij}(u_{\varepsilon\delta})}{(D_{II}(u_{\varepsilon\delta}))^{\frac{1}{2}}%
}+2\mu\varepsilon^{2} \delta^{2} D_{ij}(u_{\varepsilon\delta}), \label{sig}%
\end{equation}
where $\delta_{ij}$ is the Kronecker symbol, $g$ and $\mu$ are real positive
constants. The constant $g$ represents the yield stress of the fluid and the
constant $\mu$ is the viscosity. Relation (\ref{sig}) represents the
constitutive law of the Bingham fluid.

We define the entries of the strain tensor, denoted $D(u_{\varepsilon\delta}%
)$, by
\begin{align*}
D_{ij}(u_{\varepsilon\delta})  &  =\dfrac{1}{2}\left(  \dfrac{\partial
u_{\varepsilon\delta,i}}{\partial x_{j}}+\dfrac{\partial u_{\varepsilon
\delta,j}}{\partial x_{i}}\right)  ,1\leq i,j\leq n,\\
D_{II}(u_{\varepsilon\delta})  &  =\dfrac{1}{2}%
{\displaystyle\sum\limits_{i,j=1}^{n}}
D_{ij}(u_{\varepsilon\delta})D_{ij}(u_{\varepsilon\delta}),\\
\sigma_{ij}^{D}  &  =g\varepsilon\delta\frac{D_{ij}}{(D_{II})^{\frac{1}{2}}%
}+2\mu\varepsilon^{2}\delta^{2}D_{ij},\\
\sigma_{II}  &  =\dfrac{1}{2}%
{\displaystyle\sum\limits_{i,j=1}^{n}}
{\sigma}_{ij}^{D}{\sigma}_{ij}^{D}.
\end{align*}
Let us note that the constitutive law (\ref{sig}) is valid only if
$D_{II}(u_{\varepsilon\delta})\neq0.$ In \cite{6} it is shown that this
constitutive law is equivalent with the following one:%
\[
\left\{
\begin{array}
[c]{rcl}%
(\sigma_{II})^{\frac{1}{2}}<g\varepsilon\delta & \Leftrightarrow &
D_{ij}(u_{\varepsilon\delta})=0\\
(\sigma_{II})^{\frac{1}{2}}\geq g\varepsilon\delta & \Leftrightarrow &
D_{ij}(u_{\varepsilon\delta})=\dfrac{1}{2\mu\varepsilon^{2}\delta^{2}}\left(  1-\dfrac{g\varepsilon
\delta}{(\sigma_{II}^{\varepsilon})^{\frac{1}{2}}}\right)  {\sigma}_{ij}^{D}.
\end{array}
\right.
\]

We see that this is a threshold law: as long as the shear stress is below
$g\varepsilon\delta$, the fluid behaves as a rigid solid. When the value of
the shear stress exceeds $g\varepsilon\delta$, the fluid flows and obeys a non
linear law.

Moreover, the fluid is incompressible, which means that its velocity is
divergence free%
\[
\operatorname{div}u_{\varepsilon\delta}=0\,\,\text{in}\,\,\Omega
_{\varepsilon\delta}.
\]

In \cite{6} it is shown that the velocity $u_{\varepsilon\delta}$ satisfies
the following variational inequality when we apply to the porous media an
external force denoted by $f$ and belonging to $\left(  L^{2}\left(
\Omega\right)  \right)  ^{n}$:%
\begin{equation}
\left\{
\begin{array}
[c]{l}%
a_{\varepsilon\delta}\left(  u_{\varepsilon\delta},v-u_{\varepsilon\delta
}\right)  +j_{\varepsilon\delta}\left(  v\right)  -j_{\varepsilon\delta
}\left(  u_{\varepsilon\delta}\right)  \geq\left(  f,v-u_{\varepsilon\delta
}\right)  _{\Omega_{\varepsilon\delta}},\ \forall v\in V\left(  \Omega
_{\varepsilon\delta}\right) \\
u_{\varepsilon\delta}\in V\left(  \Omega_{\varepsilon\delta}\right)  ,
\end{array}
\right.  \label{probl1}%
\end{equation}
where
\begin{align*}
a_{\varepsilon\delta}\left(  u,v\right)   &  =2\mu\varepsilon^{2}\delta
^{2}\int_{\Omega_{\varepsilon\delta}}\nabla u\cdot\nabla
vdx,\ \ j_{\varepsilon\delta}\left(  v\right)  =g\varepsilon\delta\int
_{\Omega_{\varepsilon\delta}}\left\vert \nabla v\right\vert dx,\ \ \left(
u,v\right)  _{\Omega_{\varepsilon\delta}}=\int_{\Omega_{\varepsilon\delta}%
}u\cdot vdx,\\
V\left(  \Omega_{\varepsilon\delta}\right)   &  =\left\{  v\in\left(
H_{0}^{1}\left(  \Omega_{\varepsilon\delta}\right)  \right)  ^{n}%
:\operatorname{div}v=0\text{ in }\Omega_{\varepsilon\delta}\right\}  .
\end{align*}
If $f\in\left(  L^{2}\left(  \Omega\right)  \right)  ^{n}$, we know from
\cite{6} that for $n=2$ or $3$ and every fixed $\varepsilon$ and $\delta$
there exists a unique $u_{\varepsilon\delta}\in V\left(  \Omega_{\varepsilon
\delta}\right)  $ solution of problem (\ref{probl1}) and that if
$p_{\varepsilon\delta}$ is the pressure of the fluid in $\Omega_{\varepsilon
\delta}$, then the problem (\ref{probl1}) is equivalent to the following one:%
\begin{equation}
\left\{
\begin{array}
[c]{l}%
a_{\varepsilon\delta}\left(  u_{\varepsilon\delta},v-u_{\varepsilon\delta
}\right)  +j_{\varepsilon\delta}\left(  v\right)  -j_{\varepsilon\delta
}\left(  u_{\varepsilon\delta}\right)  \geq\left(  f,v-u_{\varepsilon\delta
}\right)  _{\Omega_{\varepsilon\delta}}-\langle\nabla p_{\varepsilon\delta},
v-u_{\varepsilon\delta} \rangle_{\Omega_{\varepsilon\delta}},\ \forall
v\in\left(  H_{0}^{1}\left(  \Omega_{\varepsilon\delta}\right)  \right)
^{n}\\
u_{\varepsilon\delta}\in V\left(  \Omega_{\varepsilon\delta}\right)
,\ p_{\varepsilon\delta}\in L_{0}^{2}\left(  \Omega_{\varepsilon\delta
}\right)
\end{array}
\right.  \label{probl2}%
\end{equation}
Here $L_{0}^{2}\left(  \Omega_{\varepsilon\delta}\right)  $ denotes the space
of functions belonging to $L^{2}\left(  \Omega_{\varepsilon\delta}\right)  $
and of mean value zero. For an open set $D$, the brackets $\langle\cdot
,\cdot\rangle_{D}$ denote the duality product between the spaces
$H^{-1}(D)^{n}$ and $H^{1}_{0}(D)^{n}$, where $H^{-1}(D)^{n}$ denotes the dual
of $H^{1}_{0}(D)^{n}$.

Our aim now is to pass to the limit as $\varepsilon\rightarrow0$ and
$\delta\rightarrow0$ in problem (\ref{probl2}). In order to do this, we first
need to get \textit{a priori} estimates for the velocity $u_{\varepsilon
\delta}$ and the pressure $p_{\varepsilon\delta}$. An important role is played
by the value of the constant in Poincar\'{e}'s inequality, with reads:

\begin{proposition}
Let $v$ be a function in $\left(  H_{0}^{1}\left(  \Omega_{\varepsilon\delta
}\right)  \right)  ^{n}$. Then we have the following inequality:
\end{proposition}

\[
\left\Vert v\right\Vert _{L^{2}\left(  \Omega_{\varepsilon\delta}\right)
^{n}}\leq C\varepsilon\delta\left\Vert \nabla v\right\Vert _{L^{2}\left(
\Omega_{\varepsilon\delta}\right)  ^{n\times n}}.
\]

\textbf{Proof.} We prove this result by using a crucial result of Tartar (see
\cite{8}), that we generalize here to the case of a domain with two scales of
periodicity. The idea is to derive Poincar\'{e}'s in the whole domain by
succesively using the $\varepsilon$-periodicity and $\delta$-periodicity of
the domain respectively and by applying the classical Poincar\'e inequality in
the cell $Z^{\ast}$.

More precisely, due to the $\varepsilon$-periodicity, it is clear that we have:%

\[
\int_{\Omega_{\varepsilon\delta}}\vert v \vert^{2} dx \approx\sum_{k
\in{\mathbb{Z}}^{n}} \int_{\varepsilon(ky^{0}+Y^{\ast}) \cap\Omega} \vert v
\vert^{2} dx
\]

and%

\[
\int_{\Omega_{\varepsilon\delta}}\vert\nabla v \vert^{2} dx \approx\sum_{k
\in{\mathbb{Z}}^{n}} \int_{\varepsilon(ky^{0}+Y^{\ast}) \cap\Omega}
\vert\nabla v \vert^{2} dx.
\]

In this above sum there are $N_{\varepsilon}$ terms and by construction
$\displaystyle N_{\varepsilon}= \frac{\vert\Omega\vert} {\vert\varepsilon Y
\vert} \approx\varepsilon^{-n} \frac{\vert\Omega\vert}{\vert Y \vert}.$

Therefore, in order to obtain the Poincar\'e inequality in the whole domain
$\Omega_{\varepsilon\delta}$ it is enough to know it in an arbitrary cell
$\varepsilon(ky^{0} +Y^{\ast})$ and then to sum over $k$.

For $v$ a function in $\left(  H_{0}^{1}\left(  \Omega_{\varepsilon\delta
}\right)  \right)  ^{n}$ we define%

\[
\displaystyle v_{\varepsilon, k^{^{\prime}}}(y)=v(\varepsilon(k^{^{\prime}%
}y^{0}+y)),\hbox { where } y \in Y^{\ast}, \, k^{^{\prime}} \in{\mathbb{Z}%
}^{n}.
\]

This function is defined on $Y^{\ast}$, it belongs to the space $H^{1}%
(Y^{\ast})$ and $v_{\varepsilon, k^{^{\prime}}}(y)=0$ on $\partial Y_{s}.$

Moreover, due to the equalities
\[
\int_{\varepsilon(hy^{0}+ Y^{\ast})} \vert v(x)\vert^{2} dx = \int_{Y^{\ast}%
}\vert v_{\varepsilon, k}(y)\vert^{2} \varepsilon^{n} dy
\]
and
\[
\int_{\varepsilon(hy^{0}+ Y^{\ast})} \vert\nabla v(x)\vert^{2} dx =
\int_{Y^{\ast}}\vert\varepsilon^{-1} \nabla_{y} v_{\varepsilon, k}(y)\vert^{2}
\varepsilon^{n} dy,
\]

it is now enough to know the Poincar\'e inequality in the domain $Y^{\ast}$ in
order to get the result.

Due to the $\delta$-periodicity and to the hypothesis on the geometry of our
domain, we have%

\[
\int_{Y^{\ast}}\vert v_{\varepsilon, k} \vert^{2} dy = \sum_{l \in{\mathbb{Z}%
}^{n}} \int_{\varepsilon(lz^{0}+Z^{\ast}) \cap Y^{\ast}} \vert v_{\varepsilon,
k} \vert^{2} dy
\]

and%

\[
\int_{Y^{\ast}}\vert\nabla_{y} v_{\varepsilon, k} \vert^{2} dy = \sum_{l
\in{\mathbb{Z}}^{n}} \int_{\varepsilon(lz^{0}+Z^{\ast}) \cap Y^{\ast}}
\vert\nabla_{y} v_{\varepsilon, k} \vert^{2} dx.
\]

In this above sum there are $N_{\delta}$ terms and by construction
$\displaystyle N_{\delta}= \frac{\vert Y^{\ast} \vert} {\vert\delta Z \vert} =
\delta^{-n} \frac{\vert Y^{\ast} \vert}{\vert Z \vert}.$

Therefore, in order to obtain the Poincar\'e inequality in the domain
$Y^{\ast}$ it is enough to know it in an arbitrary cell $\delta(lz^{0}
+Z^{\ast})$ and then to sum over $l$. By using an argument as above it is
actually enough to know the Poincar\'e inequality in the domain $Z^{\ast}.$ We define%

\[
\displaystyle v_{\varepsilon\delta, kl}(z)= v_{\varepsilon, k}(\delta
(l^{^{\prime}}z^{0}+z)),\hbox { where } z \in Z^{\ast}, \,l^{^{\prime}}
\in{\mathbb{Z}}^{n}.
\]

In $Z^{\ast}$ we know the classical Poincar\'e inequality:%
\[
\int_{Z^{\ast}}\vert v_{\varepsilon\delta, kl} \vert^{2}dz\leq
C(Z^{\ast})\int_{Z^{\ast}}\vert\nabla_{z}v_{\varepsilon\delta, kl}
\vert^{2}dz,\,\forall v\in H^{1}(Z^{\ast}), \ \ v_{\varepsilon\delta, kl}
=0\hbox { on }\partial Z_{s}.
\]

We point out that the constant $C(Z^{*})$ is independent on $\varepsilon$ and
on $\delta$. This implies%

\[
\int_{\delta(lz^{0}+Z^{\ast})}\vert v_{\varepsilon, k}(y) \vert^{2}dy=
\int_{Z^{\ast}}\vert v_{\varepsilon\delta, kl}(z) \vert^{2} \delta^{n}
dz\leq C(Z^{\ast})\delta^{2} \int_{Z^{\ast}}\vert\delta^{-1}\nabla
_{z}v_{\varepsilon\delta, kl}(z) \vert^{2} \delta^{n} dz
\]

\[
=C(Z^{\ast})\delta^{2} \int_{\delta(lz^{0}+Z^{\ast})}\vert\nabla
_{y}v_{\varepsilon, k}(y) \vert^{2} dy.
\]

By summing now over $l$ and then by repeating the same argument at the scale
$\varepsilon$ and summing over $k$ we obtain the desired result.
$\blacksquare$

\bigskip

\begin{proposition}
\label{prop3} The solution $(u_{\varepsilon\delta}, p_{\varepsilon\delta})$ of
problem (\ref{probl2}) satisfies the following \textit{a priori} estimates:%

\[
\left\Vert u_{\varepsilon\delta}\right\Vert _{L^{2}\left(  \Omega
_{\varepsilon\delta}\right)  ^{n}} \leq C
\]

\[
\varepsilon\delta\left\Vert \nabla u_{\varepsilon\delta}\right\Vert
_{L^{2}\left(  \Omega_{\varepsilon\delta}\right)  ^{n\times n}} \leq C
\]

\[
\left\Vert \nabla p_{\varepsilon\delta}\right\Vert _{H^{-1}\left(
\Omega_{\varepsilon\delta}\right)  ^{n}}\leq C\varepsilon\delta
\]

\[
\vert\vert p_{\varepsilon\delta} \vert\vert_{{L_{0}^{2}(\Omega_{\varepsilon}%
)}} \, \leq\, C.
\]

\end{proposition}

\textbf{Proof.} Setting $v=2u_{\varepsilon\delta}$ and $v=0$ successively in
(\ref{probl1}) and using the Poincar\'{e} inequality, we find the first two
estimates, for the velocity.

Let $v_{\varepsilon\delta}\in\left(  H_{0}^{1}\left(  \Omega_{\varepsilon
\delta}\right)  \right)  ^{n}$. Setting $v=v_{\varepsilon\delta}%
+u_{\varepsilon\delta}$ in (\ref{probl2}) and using estimates on the velocity,
we obtain the first estimate for the pressure and then we deduce the second
one, by using a rescaled Ne\v{c}as inequality. $\blacksquare$

Now we extend the velocity $u_{\varepsilon\delta}$ by zero to $\Omega
\diagdown\Omega_{\varepsilon\delta}$, denote the extension by the same symbol
and we have the following estimates:%
\begin{align*}
\left\Vert u_{\varepsilon\delta}\right\Vert _{L^{2}\left(  \Omega\right)
^{n}}  &  \leq C,\\
\varepsilon\delta\left\Vert \nabla u_{\varepsilon\delta}\right\Vert
_{L^{2}\left(  \Omega\right)  ^{n\times n}}  &  \leq C.
\end{align*}
Moreover, we remark that $\operatorname{div}u_{\varepsilon\delta}=0$ in
$\Omega$.

\bigskip

In order to define the extension of the pressure to the whole domain $\Omega$,
we generalize here the results from R. Bunoiu, J. Saint Jean Paulin
\cite{BuSJP}, which followed the classical idea of L. Tartar \cite{8}. We
first construct a restriction operator $S_{\varepsilon\delta}$ from
$(H_{0}^{1}\left(  \Omega\right)  )^{n}$ to $(H_{0}^{1}\left(  \Omega
_{\varepsilon\delta}\right)  )^{n}$ and using this operator we then define an
extension for the pressure to the whole domain $\Omega$.

We define the spaces $H_{s}^{1}\left(  Y^{\ast}\right)  $ and $H_{s}^{1}
\left(  Y_{f} \right)  $ by%

\[
H_{s}^{1}\left(  Y^{\ast}\right)  =\left\{  \phi\in H^{1}\left(  Y^{\ast
}\right)  :\phi=0\text{ on }\partial Y_{s}\right\}
\]

and
\[
H_{s}^{1}\left(  Y_{f}\right)  =\left\{  \phi\in H^{1}\left(  Y_{f}\right)
:\phi=0\text{ on } \Gamma\right\}  ,
\]
where the domain $Y_{f}$ is defined in section \ref{section2} and an example
is given in Figure 2. We denote $\Gamma$ the union of the boundaries of all
the obstacles contained in $Y$.

Now we first construct a restriction operator $R$ from the space
$(H^{1}(Y))^{n}$ into the space $(H_{s}^{1}\left(  Y^{\ast}\right)  )^{n} $
and next we construct a second restriction operator $W_{\delta}$ from the
space $(H_{s}^{1}\left(  Y^{\ast}\right)  )^{n} $ into the space $(H_{s}%
^{1}\left(  Y_{f}\right)  )^{n}$. Using the operators $R$ and $W_{\delta}$, we
then construct the operator%
\[
S_{\delta}: (H^{1}(Y))^{n}\rightarrow(H^{1}\left(  Y_{f}\right)  )^{n}
\]
and finally we define $S_{\varepsilon\delta}$ by applying $S_{\delta}$ to each
period $\varepsilon Y$ of $\Omega$. So we construct $S_{\varepsilon\delta}$ in
three steps, corresponding to the three following lemmas.

\begin{lemma}
\label{lemma3.1} There exists a restriction operator
\[
R:(H^{1}(Y))^{n}\rightarrow(H_{s}^{1}\left(  Y^{\ast}\right)  )^{n}%
\]
such that for $v\in(H^{1}(Y))^{n}$ we have
\end{lemma}

\begin{enumerate}
\item $Rv=v$ if $v=0$ in $Y_{s};$

\item $\operatorname{div}Rv=0$ in $Y^{\ast}$ if $\operatorname{div}v=0$ in
$Y,$

\item $\left\Vert Rv\right\Vert _{(H_{s}^{1}(Y^{\ast}))^{n}}\leq c\left\Vert
v\right\Vert _{(H^{1}(Y))^{n}}.$
\end{enumerate}

\begin{lemma}
\label{lemma3.2}There exists a restriction operator
\[
W_{\delta}:(H_{s}^{1}(Y^{\ast}))^{n}\rightarrow(H_{s}^{1}\left(  Y_{f}\right)
)^{n}%
\]
such that for $Rv\in(H_{s}^{1}(Y^{\ast}))^{n}$ we have
\end{lemma}

\begin{enumerate}
\item $W_{\delta}\left(  Rv\right)  =Rv$ if $Rv=0$ in $\bigcup\limits_{l\in
{\mathbb{Z}}^{n}}({\delta(\varepsilon)(lz^{0}+{Z}_{s})}\bigcap Y^{*}$.

\item $\operatorname{div}W_{\delta}\left(  Rv\right)  =0$ in $Y_{f}$ if
$\operatorname{div}Rv=0$ in $Y^{\ast},$

\item $\delta\left\Vert \nabla W_{\delta}\left(  Rv\right)  \right\Vert
_{\left(  L^{2}(Y_{f})\right)  ^{n\times n}}+c\left\Vert W_{\delta}\left(
Rv\right)  \right\Vert _{(L^{2}(Y_{f}))^{n}}\leq c\left\Vert v\right\Vert
_{(H^{1}(Y))^{n}}.$
\end{enumerate}

\begin{lemma}
\label{lemma3.3}There exists a restriction operator
\[
S_{\varepsilon\delta}:(H_{0}^{1}(\Omega))^{n}\rightarrow(H_{0}^{1}%
(\Omega_{\varepsilon\delta})^{n}%
\]
such that
\end{lemma}

\begin{enumerate}
\item $S_{\varepsilon\delta}\left(  v\right)  =v$ in $\Omega_{\varepsilon
\delta},$ $\forall v\in(H_{0}^{1}(\Omega_{\varepsilon\delta}))^{n},$

\item $\operatorname{div}S_{\varepsilon\delta}v=0$ in $\Omega_{\varepsilon
\delta}$ if $\operatorname{div}v=0$ in $\Omega,$

\item $\displaystyle  \varepsilon\delta\left\Vert \nabla S_{\varepsilon\delta
}v\right\Vert _{\left(  L^{2}(\Omega_{\varepsilon\delta})\right)  ^{n \times
n}} + c \left\Vert S_{\varepsilon\delta}v\right\Vert _{(L^{2}(\Omega
_{\varepsilon\delta}))^{n}} \leq c \left\Vert v \right\Vert _{(H^{1}%
(\Omega))^{n}}. $
\end{enumerate}

Let now $v$ be a function in the space $H_{0}^{1}(\Omega).$ As $\nabla
p_{\varepsilon\delta}\in H^{-1}\left(  \Omega_{\varepsilon\delta}\right)  $,
we define the application $F^{\varepsilon}$ by%
\[
\left\langle F^{\varepsilon},v\right\rangle _{\Omega}=\left\langle \nabla
p_{\varepsilon\delta},S_{\varepsilon\delta}v\right\rangle _{\Omega
_{\varepsilon\delta}},
\]
where $S_{\varepsilon\delta}$ is the operator defined by Lemma \ref{lemma3.3}.
The following proposition defines us the extension $\widetilde{p}%
_{\varepsilon\delta}$ of the pressure $p_{\varepsilon\delta}$ to the whole
$\Omega$. Moreover, it gives us a strong convergence result for this
extension. Following the ideas of L.Tartar \cite{8}, we can prove

\begin{proposition}
Let $p_{\varepsilon\delta}$ be as in (\ref{probl2}). Then, for each
$\varepsilon$ and $\delta$ there exists an extension $\widetilde
{p}_{\varepsilon\delta}$ of $p_{\varepsilon\delta}$ defined on $\Omega$ such
that%
\[
\widetilde{p}_{\varepsilon\delta}=p_{\varepsilon\delta}\text{ in }%
\Omega_{\varepsilon\delta}.
\]
Moreover, up to a subsequence, we have%
\[
\widetilde{p}_{\varepsilon\delta}\rightarrow\widehat{p}\text{ strongly in
}L^{2}_{0}\left(  \Omega\right)  .
\]
The function $F^{\varepsilon}$ defined before and the pressure $\widetilde
{p}_{\varepsilon\delta}$ are linked by the relation
\[
F^{\varepsilon}=\nabla\widetilde{p}_{\varepsilon\delta}.
\]

\end{proposition}

For every function $v$ that is the extension by zero to the whole $\Omega$ of
a function in $H_{0}^{1}(\Omega_{\varepsilon\delta})^{n}$ we deduce:%

\[
-\langle\nabla p_{\varepsilon\delta},v\rangle_{\Omega_{\varepsilon\delta}%
}=-\langle\nabla\widetilde{p}_{\varepsilon\delta},v\rangle_{\Omega}= \left(
\widetilde{p}_{\varepsilon\delta},\operatorname{div}v\right)  _{\Omega}.
\]

According to the extensions of the velocity and of the pressure, problem
(\ref{probl2}) can now be written as:
\begin{align}
&  2\mu\varepsilon^{2}\delta^{2}\int_{\Omega}\nabla u_{\varepsilon\delta}%
\cdot\nabla\left(  v-u_{\varepsilon\delta}\right)  dx+g\varepsilon\delta
\int_{\Omega}\left\vert \nabla v\right\vert dx-g\varepsilon\delta\int_{\Omega
}\left\vert \nabla u_{\varepsilon\delta}\right\vert dx\label{ineq}\\
&  \geq\int_{\Omega}f_{\varepsilon}\left(  v-u_{\varepsilon\delta}\right)
dx+\int_{\Omega}\widetilde{p}_{\varepsilon\delta}\operatorname{div}\left(
v-u_{\varepsilon\delta}\right)  dx,\nonumber
\end{align}
for every $v$ that is the extension by zero to the whole $\Omega$ of a
function in $H_{0}^{1}(\Omega_{\varepsilon\delta})^{n}$.

\section{Convergence result\label{section4}}

Now we state the main result of our paper, the convergence result for the
variational inequality (\ref{ineq}). In order to prove it, we apply the
unfolding operators from section \ref{section2}, together with the \textit{a
priori} estimates from Proposition \ref{prop3} and the compactness results
from Proposition \ref{prop2}.

\begin{theorem}
\label{theor1}Let $u_{\varepsilon\delta}$ and $\widetilde{p}_{\varepsilon
\delta}$ verify relation (\ref{ineq}). Then there exist $\widehat{u}\left(
x,y,z\right)  \in L^{2}\left(  \Omega\times Y^{\ast};(H_{per}^{1}\left(
Z^{\ast}\right)  )^{n}\right)  $ and $\widehat{p}\in L_{0}^{2}(\Omega)\cap
H^{1}(\Omega)$ such that
\[
u_{\varepsilon\delta}\rightarrow\frac{1}{\vert Y \vert\vert Z \vert}
{\displaystyle\int_{Y^{\ast}}}
{\displaystyle\int_{Z^{\ast}}}
\widehat{u}\left(  \cdot,y,z\right)  dydz \hbox { weakly in } (L^{2}\left(
\Omega\right)  )^{n},
\]

\[
\widetilde{p}_{\varepsilon\delta}\rightarrow\widehat{p} \hbox { strongly in }
L_{0}^{2}(\Omega)
\]

and satisfy the limit problem
\begin{align}
&  2\mu%
{\displaystyle\int_{Y^{\ast}}}
{\displaystyle\int_{Z^{\ast}}}
\nabla_{z}\widehat{u}\cdot\nabla_{z}\left(  \Phi-\widehat{u}\right)  dydz+g%
{\displaystyle\int_{Y^{\ast}}}
{\displaystyle\int_{Z^{\ast}}}
\left\vert \nabla_{z}\left(  \Phi\right)  \right\vert dydz-g%
{\displaystyle\int_{Y^{\ast}}}
{\displaystyle\int_{Z^{\ast}}}
\left\vert \nabla_{z}\widehat{u}\right\vert dydz\label{limprob}\\
&  \geq\langle f-\nabla_{x}\widehat{p},\Phi-\widehat{u}\rangle_{Y^{\ast}\times
Z^{\ast}}\nonumber
\end{align}
for every $\Phi\in\mathcal{V},$ where%
\[
\mathcal{V}= \left\{  \Phi\in L^{2}(Y^{\ast};\mathcal{W}),\text{
}\operatorname{div}_{y}\int_{Z^{\ast}}\Phi dz=0\,\hbox {in }\,Y^{\ast
},\text{ }\int_{Z^{\ast}}\Phi dz\,\hbox {is Y- periodic},\text{ }%
\nu_{Y}\cdot\int_{Z^{\ast}}\Phi dz=0\,\hbox {on}\,\partial
Y_{s}\right\}  .
\]%
\[
\mathcal{W}=\{\phi\in H_{per,0}^{1}(Z^{\ast}),\text{ }\operatorname{div}%
_{z}\phi=0\,\hbox {in }\,Z^{\ast}\}
\]%
\[
H_{per,0}^{1}(Z^{\ast})=\{\phi\in(H^{1}(Z^{\ast}))^{n},\text{ }\phi
=0\,\hbox { on }\partial Z_{s},\text{ }\phi\,\hbox { is Z periodic}\}.
\]
The function $\widehat{u}$ satisfies the following conditions:%
\begin{align}
\operatorname{div}_{x}%
{\displaystyle\int_{Y^{\ast}}}
{\displaystyle\int_{Z^{\ast}}}
\widehat{u}dydz  &  =0\text{ in }\Omega,\label{5.2}\\
\operatorname{div}_{y}%
{\displaystyle\int_{Z^{\ast}}}
\widehat{u}dz  &  =0\text{ in }\Omega\times Y^{\ast},\label{5.3}\\
\operatorname{div}_{z}\widehat{u}  &  =0\text{ in }\Omega\times Y^{\ast}\times
Z^{\ast}\label{5.4}\\
\nu\cdot%
{\displaystyle\int_{Y^{\ast}}}
{\displaystyle\int_{Z^{\ast}}}
\widehat{u}dydz  &  =0\text{ on }\partial\Omega,\label{5.5}\\
&  \widehat{u}\cdot\nu_{Z}\text{ takes opposite values on opposite faces of
}Z\label{5.6}\\
\nu_{Y}\cdot%
{\displaystyle\int_{Z^{\ast}}}
\widehat{u}dz  &  =0\text{ takes opposite values on opposite faces of
}Y,\label{5.7}\\
\widehat{u}\cdot\nu_{Z}  &  =0\text{ on }\partial Z_{s}\label{5.8}\\
\nu_{Y}\cdot%
{\displaystyle\int_{Z^{\ast}}}
\widehat{u}dz  &  =0\text{ on }\partial Y_{s} \label{5.9}%
\end{align}

\end{theorem}

\textbf{Proof.} Taking into account the \textit{a priori} estimates from
Proposition \ref{prop3} and then using Proposition \ref{prop2} and Remark
\ref{rk}, we have the following convergences for the velocity:
\begin{align*}
\left\Vert u_{\varepsilon\delta}\right\Vert _{L^{2}\left(  \Omega\right)
^{n}}  &  \leq C\Rightarrow T_{\delta}\left(  T_{\varepsilon}\left(
u_{\varepsilon\delta}\right)  \right)  \rightarrow\widehat{u}\text{ weakly in
}(L^{2}\left(  \Omega\times Y^{\ast}\times Z^{\ast}\right)  )^{n},\\
\varepsilon\delta\left\Vert \nabla u_{\varepsilon\delta}\right\Vert
_{L^{2}\left(  \Omega\right)  ^{n\times n}}  &  \leq C\Rightarrow
\varepsilon\delta T_{\delta}\left(  T_{\varepsilon}\left(  u_{\varepsilon
\delta}\right)  \right)  \rightarrow\nabla_{z}\widehat{u}\text{ weakly in
}(L^{2}\left(  \Omega\times Y^{\ast}\times Z^{\ast}\right)  )^{n\times n},\\
\widehat{u}  &  \in L^{2}\left(  \Omega\times Y^{\ast};(H_{per}^{1}\left(
Z^{\ast}\right)  )^{n}\right)  ,\\
u_{\varepsilon\delta}  &  \rightarrow\int_{Y^{\ast}\times Z^{\ast}}\widehat
{u}\left(  \cdot,y,z\right)  dydz\text{ weakly in }(L^{2}\left(
\Omega\right)  )^{n}.
\end{align*}
According to \cite{8}, we have for the pressure the converegnce
\[
\widetilde{p}_{\varepsilon\delta}\rightarrow\widehat{p}\text{ strongly in
}L_{0}^{2}(\Omega).
\]
Using property $p_{4}$) of the unfolding operators we get:%
\[
T_{\delta}\left(  T_{\varepsilon}\left(  \widetilde{p}_{\varepsilon\delta
}\right)  \right)  \rightarrow\widehat{p}\text{ strongly in }L_{0}^{2}\left(
\Omega\times Y \times Z \right)  .
\]

In order to prove relation (\ref{5.2}), let us observe that
$\operatorname{div}u_{\varepsilon\delta}=0$ implies $\varepsilon
T_{\varepsilon}\left(  \operatorname{div}u_{\varepsilon\delta}\right)  =0.$
But
\[
\varepsilon T_{\varepsilon}\left(  \operatorname{div}u_{\varepsilon\delta
}\right)  =\varepsilon T_{\varepsilon}\left(  \sum_{i=1}^{n}\frac{\partial
u_{\varepsilon\delta,i}}{\partial x_{i}}\right)  =\varepsilon T_{\varepsilon
}\left(  \sum_{i=1}^{n}\frac{1}{\varepsilon}\frac{\partial u_{\varepsilon
\delta,i}}{\partial y_{i}}\right)  =\operatorname{div}_{y}T_{\varepsilon
}\left(  u_{\varepsilon\delta}\right)
\]
and so
\[
\delta\varepsilon T_{\delta}(T_{\varepsilon}\left(  \operatorname{div}%
u_{\varepsilon\delta}\right)  )=\delta T_{\delta}(\varepsilon T_{\varepsilon
}(\operatorname{div}u_{\varepsilon\delta}))=\delta T_{\delta}%
(\operatorname{div}_{y}T_{\varepsilon}(u_{\varepsilon\delta}%
)=\operatorname{div}_{z}T_{\delta}(T_{\varepsilon}\left(  u_{\varepsilon
\delta}\right)  )
\]
which implies $\operatorname{div}_{z}T_{\delta}(T_{\varepsilon}\left(
u_{\varepsilon\delta}\right)  )=0.$

We pass to the limit as $\varepsilon$ tends to zero in this last equality and
we get the desired result.

In order to prove (\ref{5.3}) let us take $\Psi\in\mathcal{D}\left(
\Omega\right)  $, $\psi\in H_{per}^{1}(Y)$ and define $\psi_{\varepsilon
}(x)=\psi\left(  \frac{x}{\varepsilon}\right)  .$

We have
\[
0=\int_{\Omega}\varepsilon\operatorname{div}u_{\varepsilon\delta}%
\Psi\psi_{\varepsilon}dx=\int_{\Omega}\varepsilon u_{\varepsilon\delta
}(\nabla_{x}\Psi\psi_{\varepsilon}+\Psi\nabla_{x}\psi_{\varepsilon})dx.
\]
By applying the unfolding at the scale $\varepsilon$ we get%

\[
0=\int_{\Omega}\int_{Y}T_{\varepsilon}(u_{\varepsilon\delta
})(\varepsilon\nabla_{x}\Psi T_{\varepsilon}(\psi_{\varepsilon})+\Psi
T_{\varepsilon}(\varepsilon\nabla_{x}\psi_{\varepsilon}))dxdy.
\]
We pass to the limit as $\varepsilon$ tends to zero and we get%

\[
0=\int_{\Omega}\int_{Y}\Bigl(\frac{1}{|Z|}\int_{Z}%
\widehat{u}(x,y,z)dz\Bigr)\nabla_{y}\psi(y)\Psi(x)dxdy.
\]
An integration by parts in the domain $Y$ gives%

\[
0=\int_{\Omega}\int_{Y}\operatorname{div}_{y}\Bigl(\frac{1}%
{|Z|}\int_{Z}\widehat{u}(x,y,z)dz\Bigr)\psi(y)\Psi(x)dxdy,
\]
and this last equality implies (\ref{5.3}).

In order to prove relation (\ref{5.4}), let us take $\Psi\in\mathcal{D}\left(
\Omega\right)  .$

We have%
\[
0=\int_{\Omega}\operatorname{div}u_{\varepsilon\delta}\Psi dx=\int_{\Omega
}u_{\varepsilon\delta}\nabla\Psi dx.
\]
By applying the unfolding at scale $\varepsilon$ and then at scale $\delta$ we
get%
\[
0=\int_{\Omega\times Y\times Z}T_{\delta}(T_{\varepsilon}\left(
u_{\varepsilon\delta}\right)  )T_{\delta}(T_{\varepsilon}\left(  \nabla
\Psi\right)  )dxdydz.
\]
We pass to the limit as $\varepsilon$ tends to zero and we get%
\begin{align*}
0  &  =\int_{\Omega\times Y\times Z}\widehat{u}\,\nabla_{x}\Psi dxdydz,\\
0  &  =\int_{\Omega}\operatorname{div}_{x}\left(  \int_{Y\times Z}\widehat
{u}\left(  x,y,z\right)  dydz\right)  \Psi dx,\ \ \ \forall\Psi\in
\mathcal{D}\left(  \Omega\right)  ,
\end{align*}
which implies (\ref{5.4}).

Relation (\ref{5.5}) is a consequence of the following assertions:%
\begin{align*}
\widehat{u}(x,y,z)  &  =0\text{ in }Y_{s}\times Z_{s},\text{ a.e. in }%
\Omega,\\
u_{\varepsilon\delta}  &  \rightarrow\frac{1}{\left\vert Y\right\vert
\left\vert Z\right\vert }\int_{Y\times Z}\widehat{u}\left(  x,y,z\right)
dydz\text{ weakly in }(L^{2}\left(  \Omega)\right)  ^{n},
\end{align*}
together with the linearity and continuity of the normal trace application
from the space $H\left(  \operatorname{div},\Omega\right)  =\left\{
\varphi\in\left(  L^{2}\left(  \Omega\right)  \right)  ^{n}:\operatorname{div}%
\varphi\in L^{2}\left(  \Omega\right)  \right\}  $ into $H^{-1/2}\left(
\partial\Omega\right)  .$

By choosing particular test functions in relations (\ref{5.2}) and (\ref{5.3})
we obtain relations (\ref{5.6}) and (\ref{5.7}) respectively.

Relation (\ref{5.8}) is a consequence of relations (\ref{5.2}) and (\ref{5.6}).

Relation (\ref{5.9}) is a consequence of relations (\ref{5.3}) and (\ref{5.7}).

By applying now the unfolding operator to the inequality (\ref{ineq}), we get%
\begin{align}
&  2\mu\varepsilon^{2}\delta^{2}\int_{\Omega\times Y^{\ast}\times Z^{\ast}%
}T_{\delta}\left(  T_{\varepsilon}\left(  \nabla u_{\varepsilon\delta}\right)
\right)  \cdot T_{\delta}\left(  T_{\varepsilon}\left(  \nabla\left(
v-u_{\varepsilon\delta}\right)  \right)  \right)  dxdydz\label{20}\\
&  +g\varepsilon\delta\int_{\Omega\times Y^{\ast}\times Z^{\ast}}T_{\delta
}\left(  T_{\varepsilon}\left(  \left\vert \nabla v\right\vert \right)
\right)  dxdydz-g\varepsilon\delta\int_{\Omega\times Y^{\ast}\times Z^{\ast}%
}T_{\delta}\left(  T_{\varepsilon}\left(  \left\vert \nabla u_{\varepsilon
\delta}\right\vert \right)  \right)  dxdydz\nonumber\\
&  \geq\int_{\Omega\times Y^{\ast}\times Z^{\ast}}T_{\delta}\left(
T_{\varepsilon}\left(  f_{\varepsilon}\right)  \right)  T_{\delta}\left(
T_{\varepsilon}\left(  v-u_{\varepsilon\delta}\right)  \right)  dxdydz+\int
_{\Omega\times Y^{\ast}\times Z^{\ast}}T_{\delta}\left(  T_{\varepsilon
}\left(  \widetilde{p}_{\varepsilon}\right)  \right)  T_{\delta}\left(
T_{\varepsilon}\left(  \operatorname{div}\left(  v-u_{\varepsilon\delta
}\right)  \right)  \right)  dxdydz.\nonumber
\end{align}

In order to pass to the limit in relation (\ref{20}), we consider a test
function $v=v_{\varepsilon\delta}$ of the form:%
\begin{equation}
v_{\varepsilon\delta}\left(  x\right)  =\Psi\left(  x\right)  \psi\left(
\frac{x}{\varepsilon}\right)  \phi\left(  \frac{x}{\varepsilon\delta}\right)
, \label{testf}%
\end{equation}
where $\Psi\in\mathcal{D}\left(  \Omega\right)  $, $\psi\in\mathcal{D}%
(Y^{\ast}),$ $\phi\in(H_{per,0}^{1}(Z^{\ast}))^{n},$ $\operatorname{div}%
_{z}\phi=0.$

We have
\begin{align}
\nabla_{x}v_{\varepsilon\delta}  &  =\nabla_{x}\left(  \Psi\left(  x\right)
\psi\left(  \frac{x}{\varepsilon}\right)  \phi\left(  \frac{x}{\varepsilon
\delta}\right)  \right)  =\nabla_{x}\Psi\left(  x\right)  \psi\left(  \frac
{x}{\varepsilon}\right)  \phi\left(  \frac{x}{\varepsilon\delta}\right)
+\frac{1}{\varepsilon}\Psi\left(  x\right)  \nabla_{x}\psi\left(  \frac
{x}{\varepsilon}\right)  \phi\left(  \frac{x}{\varepsilon\delta}\right)
\label{testf2}\\
&  +\frac{1}{\varepsilon\delta}\Psi\left(  x\right)  \psi\left(  \frac
{x}{\varepsilon}\right)  \nabla_{x}\phi\left(  \frac{x}{\varepsilon\delta
}\right)  .\nonumber
\end{align}
By using this test function we get for the first term in relation (\ref{20}):%
\begin{align*}
&  2\mu\varepsilon^{2}\delta^{2}\int_{\Omega\times Y^{\ast}\times Z^{\ast}%
}T_{\delta}\left(  T_{\varepsilon}\left(  \nabla u_{\varepsilon\delta}\right)
\right)  \cdot T_{\delta}\left(  T_{\varepsilon}\left(  \nabla\left(
v_{\varepsilon\delta}-u_{\varepsilon\delta}\right)  \right)  \right)  dxdydz\\
&  =2\mu\varepsilon^{2}\delta^{2}\int_{\Omega\times Y^{\ast}\times Z^{\ast}%
}T_{\delta}\left(  T_{\varepsilon}\left(  \nabla u_{\varepsilon\delta}\right)
\right)  \cdot T_{\delta}\left(  T_{\varepsilon}\left(  \nabla\left(
v_{\varepsilon\delta}\right)  \right)  \right)  dxdydz\\
&  -2\mu\varepsilon^{2}\delta^{2}\int_{\Omega\times Y^{\ast}\times Z^{\ast}%
}T_{\delta}\left(  T_{\varepsilon}\left(  \nabla u_{\varepsilon\delta}\right)
\right)  \cdot T_{\delta}\left(  T_{\varepsilon}\left(  \nabla\left(
u_{\varepsilon\delta}\right)  \right)  \right)  dxdydz\\
&  =2\mu\varepsilon^{2}\delta^{2}\int_{\Omega\times Y^{\ast}\times Z^{\ast}%
}T_{\delta}\left(  T_{\varepsilon}\left(  \nabla u_{\varepsilon\delta}\right)
\right)  \cdot\Big[T_{\delta}\left(  T_{\varepsilon}\left(  \nabla_{x}%
\Psi\right)  \right)  T_{\delta}\left(  T_{\varepsilon}\left(  \psi\right)
\right)  T_{\delta}\left(  T_{\varepsilon}\left(  \phi\right)  \right) \\
&  \left.  +\frac{1}{\varepsilon}T_{\delta}\left(  T_{\varepsilon}\left(
\Psi\right)  \right)  T_{\delta}\left(  T_{\varepsilon}\left(  \nabla_{y}%
\psi\right)  \right)  T_{\delta}\left(  T_{\varepsilon}\left(  \phi\right)
\right)  +\frac{1}{\varepsilon\delta}T_{\delta}\left(  T_{\varepsilon}\left(
\Psi\right)  \right)  T_{\delta}\left(  T_{\varepsilon}\left(  \psi\right)
\right)  T_{\delta}\left(  T_{\varepsilon}\left(  \nabla_{z}\phi\right)
\right)  \right]  dxdydz\\
&  -2\mu\varepsilon^{2}\delta^{2}\int_{\Omega\times Y^{\ast}\times Z^{\ast}%
}T_{\delta}\left(  T_{\varepsilon}\left(  \nabla u_{\varepsilon\delta}\right)
\right)  \cdot T_{\delta}\left(  T_{\varepsilon}\left(  \nabla u_{\varepsilon
\delta}\right)  \right)  dxdydz\\
&  =2\mu\int_{\Omega\times Y^{\ast}\times Z^{\ast}}\varepsilon\delta
T_{\delta}\left(  T_{\varepsilon}\left(  \nabla u_{\varepsilon\delta}\right)
\right)  \cdot\varepsilon\delta T_{\delta}\left(  T_{\varepsilon}\left(
\nabla_{x}\Psi\right)  \right)  \psi\left(  y\right)  \phi\left(  z\right)
dxdydz\\
&  +2\mu\int_{\Omega\times Y^{\ast}\times Z^{\ast}}\varepsilon\delta
T_{\delta}\left(  T_{\varepsilon}\left(  \nabla u_{\varepsilon\delta}\right)
\right)  \cdot\delta T_{\delta}\left(  T_{\varepsilon}\left(  \Psi\right)
\right)  \nabla_{y}\psi\left(  y\right)  \phi\left(  z\right)  dxdydz\\
&  +2\mu\int_{\Omega\times Y^{\ast}\times Z^{\ast}}\varepsilon\delta
T_{\delta}\left(  T_{\varepsilon}\left(  \nabla u_{\varepsilon\delta}\right)
\right)  \cdot T_{\delta}\left(  T_{\varepsilon}\left(  \Psi\right)  \right)
\psi\left(  y\right)  \nabla_{z}\phi\left(  z\right)  dxdydz\\
&  -2\mu\int_{\Omega\times Y^{\ast}\times Z^{\ast}}\left\vert \varepsilon
\delta T_{\delta}\left(  T_{\varepsilon}\left(  \nabla u_{\varepsilon\delta
}\right)  \right)  \right\vert ^{2}dxdydz.
\end{align*}
According to the general convergence results for the unfolding we have that
the first and second terms tend to zero and the third one to the following
limit:%
\[
2\mu\int_{\Omega\times Y^{\ast}\times Z^{\ast}}\nabla_{z}\widehat{u}\cdot
\Psi\left(  x\right)  \psi\left(  y\right)  \nabla_{z}\phi\left(  z\right)
dxdydz.
\]
By using now the fact that the function $B(\varphi)=|\varphi|^{2}$ is proper
convex continuous, we have for the fourth term
\[
\liminf_{\varepsilon\rightarrow0}2\mu\int_{\Omega\times Y^{\ast}\times
Z^{\ast}}|\varepsilon\delta T_{\delta}\left(  T_{\varepsilon}\left(  \nabla
u_{\varepsilon\delta}\right)  \right)  |^{2}dxdydz\,\geq2\mu\int_{\Omega\times
Y^{\ast}\times Z^{\ast}}|\nabla_{z}\widehat{u}|^{2}dxdydz.
\]
In order to pass to the limit in the non linear terms, let us first remark
that for a function $v$ in $\left(  H^{1}(\Omega)\right)  ^{n}$ we have
\begin{equation}
\varepsilon\delta T_{\delta}\left(  T_{\varepsilon}(|\nabla v|)\right)
=|\nabla_{z}T_{\delta}\left(  T_{\varepsilon}(v)\right)  |. \label{grad}%
\end{equation}

Indeed, according to a result in \cite{BuCaPe}, we know that%
\[
\varepsilon T_{\varepsilon}(|\nabla v|)=|\nabla_{y}T_{\varepsilon}(v)|
\]
and following the same ideas we can prove that for a function $w\in
L^{2}(\Omega,H^{1}(Y))$ we have%
\[
\delta T_{\delta}\left(  (|\nabla_{y}w|)\right)  =|\nabla_{z}T_{\delta}\left(
w\right)  |.
\]
This implies%
\[
\varepsilon\delta T_{\delta}\left(  T_{\varepsilon}(|\nabla v|)\right)
=\delta T_{\delta}(\varepsilon T_{\varepsilon}(|\nabla v|))=\delta T_{\delta
}(|\nabla_{y}T_{\varepsilon}(v)|)=|\nabla_{z}T_{\delta}\left(  T_{\varepsilon
}(v)\right)  |.
\]

In order to pass to the limit in the first non linear term, by using the
previous identity for the function $v_{\varepsilon\delta}$ given by
(\ref{testf}), we have%
\begin{align*}
&  \left\vert g\varepsilon\delta\int_{\Omega\times Y^{\ast}\times Z^{\ast}%
}T_{\delta}(T_{\varepsilon}\left(  \left\vert \nabla v_{\varepsilon\delta
}\right\vert \right)  )dxdydz-g\int_{\Omega\times Y^{\ast}\times Z^{\ast}%
}\left\vert \nabla_{y}\left(  \Psi\psi\phi\right)  \right\vert
dxdydz\right\vert \\
&  =\left\vert g\int_{\Omega\times Y^{\ast}\times Z^{\ast}}\left\vert
\nabla_{z}T_{\delta}(T_{\varepsilon}\left(  v_{\varepsilon\delta}\right)
)\right\vert dxdydz-g\int_{\Omega\times Y^{\ast}\times Z^{\ast}}\left\vert
\nabla_{z}\left(  \Psi\psi\phi\right)  \right\vert dxdydz\right\vert \\
&  \leq g\int_{\Omega\times Y^{\ast}\times Z^{\ast}}\left\vert \nabla
_{z}T_{\delta}(T_{\varepsilon}\left(  v_{\varepsilon\delta}\right)
)-\nabla_{z}\left(  \Psi\psi\phi\right)  \right\vert dxdydz\\
&  =g\int_{\Omega\times Y^{\ast}\times Z^{\ast}}|\varepsilon\delta T_{\delta
}(T_{\varepsilon}\left(  \nabla_{x}\Psi\right)  )\left(  x,y,z\right)
\psi(y)\cdot\phi\left(  z\right)  +\delta T_{\delta}(T_{\varepsilon}\left(
\Psi\right)  )\left(  x,y,z\right)  \nabla_{y}(\psi(y))\phi(z)\\
&  \ \ \ \ \ \ \ \ \ \ \ \ \ \ \ \ \ \ \ \ \ \ +T_{\delta}(T_{\varepsilon
}\left(  \Psi\right)  )\left(  x,y,z\right)  \psi(y)\nabla_{z}\phi\left(
z\right)  -\Psi\left(  x\right)  \psi(y)\nabla_{z}\phi\left(  z\right)
|dxdydz\\
&  \leq g\int_{\Omega\times Y^{\ast}\times Z^{\ast}}\left\vert T_{\delta
}(T_{\varepsilon}\left(  \varepsilon\delta\nabla_{x}\Psi\right)  \left(
x,y,z\right)  \psi(y)\cdot\phi\left(  z\right)  \right\vert dxdydz\\
&  +g\int_{\Omega\times Y^{\ast}\times Z^{\ast}}\left\vert T_{\delta}(\delta
T_{\varepsilon}\left(  \Psi\right)  \left(  x,y,z\right)  \nabla_{y}%
\psi(y)\cdot\phi\left(  z\right)  \right\vert dxdydz\\
&  +g\int_{\Omega\times Y^{\ast}\times Z^{\ast}}\left\vert \left(  T_{\delta
}(T_{\varepsilon}\left(  \Psi\right)  \left(  x,y,z\right)  -\Psi\left(
x\right)  \right)  \psi(y)\nabla_{z}\phi\left(  z\right)  \right\vert dxdydz\\
&  \leq g\left\Vert T_{\delta}(T_{\varepsilon}\left(  \varepsilon\delta
\nabla_{x}\Psi\right)  \right\Vert _{\left(  L^{2}\left(  \Omega\times
Y^{\ast}\times Z^{\ast}\right)  \right)  ^{n}}\left\Vert \psi\phi\right\Vert
_{\left(  L^{2}\left(  \Omega\times Y^{\ast}\times Z^{\ast}\right)  \right)
^{n}}\\
&  +g\left\Vert T_{\delta}(\delta T_{\varepsilon}(\Psi))\right\Vert
_{L^{2}\left(  \Omega\times Y^{\ast}\times Z^{\ast}\right)  }\left\Vert
\nabla_{y}\psi(y)\cdot\phi(z)\right\Vert _{L^{2}\left(  \Omega\times Y^{\ast
}\times Z^{\ast}\right)  }\\
&  +\left\Vert T_{\delta}(T_{\varepsilon}\left(  \Psi\right)  )-\Psi
\right\Vert _{L^{2}\left(  \Omega\times Y^{\ast}\times Z^{\ast}\right)
}\left\Vert \psi(y)\nabla_{z}\left(  \phi\right)  \right\Vert _{\left(
L^{2}\left(  \Omega\times Y^{\ast}\times Z^{\ast}\right)  \right)  ^{n\times
n}}.
\end{align*}
Passing to the limit as $\varepsilon\rightarrow0,$ we have that
\[
\left\Vert T_{\delta}(T_{\varepsilon}\left(  \varepsilon\delta\nabla_{x}%
\Psi\right)  )\right\Vert _{\left(  L^{2}\left(  \Omega\times Y^{\ast}\times
Z^{\ast}\right)  \right)  ^{n}}\rightarrow0
\]
and
\[
\left\Vert T_{\delta}(\delta T_{\varepsilon}\left(  \Psi\right)  )\right\Vert
_{L^{2}\left(  \Omega\times Y^{\ast}\times Z^{\ast}\right)  }\rightarrow0.
\]
Moreover, $T_{\delta}(T_{\varepsilon}(\Psi))\rightarrow\Psi$\ strongly in
$L^{2}\left(  \Omega\times Y^\ast\times Z^{\ast}\right)  $ and so
\[
\left\Vert T_{\delta}(T_{\varepsilon}\left(  \Psi\right)  )-\Psi\right\Vert
_{L^{2}\left(  \Omega\times Y^{\ast}\times Z^{\ast}\right)  }\rightarrow0.
\]

Then%
\[
\lim_{\varepsilon\rightarrow0}g\varepsilon\delta\int_{\Omega\times Y^{\ast
}\times Z^{\ast}}T_{\delta}\left(  T_{\varepsilon}\left(  \left\vert \nabla
v_{\varepsilon}\right\vert \right)  \right)  dxdydz=g\int_{\Omega\times
Y^{\ast}\times Z^{\ast}}\left\vert \nabla_{z}(\Psi\left(  x\right)
\psi\left(  y\right)  \phi\left(  z\right)  ) \right\vert dxdydz.
\]

In order to pass to the limit in the second non linear term, we use identity
(\ref{grad}) for the function $u_{\varepsilon\delta}$ and the fact that the
function $E(\varphi)=|\varphi|$ is proper convex continuous. We then deduce:%
\[
\liminf_{\varepsilon\rightarrow0}g\varepsilon\delta\int_{\Omega\times Y^{\ast
}\times Z^{\ast}}T_{\delta}\left(  T_{\varepsilon}\left(  \left\vert \nabla
u_{\varepsilon\delta}\right\vert \right)  \right)  dxdydz\geq g\int
_{\Omega\times Y^{\ast}\times Z^{\ast}}|\nabla_{z}\widehat{u}|dxdydz.
\]

Moreover,
\begin{align*}
&  \int_{\Omega\times Y^{\ast}\times Z^{\ast}}T_{\delta}\left(  T_{\varepsilon
}\left(  f_{\varepsilon}\right)  \right)  T_{\delta}\left(  T_{\varepsilon
}\left(  v\right)  \right)  dxdydz-\int_{\Omega\times Y^{\ast}\times Z^{\ast}%
}T_{\delta}\left(  T_{\varepsilon}\left(  f_{\varepsilon}\right)  \right)
T_{\delta}\left(  T_{\varepsilon}\left(  u_{\varepsilon\delta}\right)
\right)  dxdydz\\
&  \rightarrow\int_{\Omega\times Y^{\ast}\times Z^{\ast}}f\Psi\psi\phi
dxdydz-\int_{\Omega\times Y^{\ast}\times Z^{\ast}}f\widehat{u}dxdydz.
\end{align*}
We consider now the term ${%
{\displaystyle\int_{\Omega\times Y^{\ast}\times Z^{\ast}}}
T_{\delta}}\left(  {T_{\varepsilon}\left(  \widetilde{p}_{\varepsilon\delta
}\right)  }\right)  T_{\delta}\left(  {T_{\varepsilon}\left(
\operatorname{div}\left(  v-u_{\varepsilon\delta}\right)  \right)  }\right)
{dxdydz}$. Using $\operatorname{div}_{x}u_{\varepsilon\delta}=0$ we obtain:
\begin{align*}
&  {\int_{\Omega\times Y^{\ast}\times Z^{\ast}}T_{\delta}}\left(
{T_{\varepsilon}\left(  \widetilde{p}_{\varepsilon\delta}\right)  }\right)
T_{\delta}\left(  {T_{\varepsilon}\left(  \operatorname{div}\left(
v-u_{\varepsilon\delta}\right)  \right)  }\right)  {dxdydz}\\
&  ={\int_{\Omega\times Y^{\ast}\times Z^{\ast}}T_{\delta}}\left(
{T_{\varepsilon}\left(  \widetilde{p}_{\varepsilon\delta}\right)  }\right)
T_{\delta}\left(  {T_{\varepsilon}\left(  \operatorname{div}_{x}v\right)
}\right)  {dxdydz}\\
&  ={\int_{\Omega\times Y^{\ast}\times Z^{\ast}}T_{\delta}}\left(
{T_{\varepsilon}\left(  \widetilde{p}_{\varepsilon\delta}\right)  }\right)
T_{\delta}\left(  T_{\varepsilon}\left(  \operatorname{div}_{x}\left(
\Psi\left(  x\right)  \psi\left(  \frac{x}{\varepsilon}\right)  \phi\left(
\frac{x}{\varepsilon\delta}\right)  \right)  \right)  \right)  dxdydz\\
&  ={\int_{\Omega\times Y^{\ast}\times Z^{\ast}}T_{\delta}}\left(
{T_{\varepsilon}\left(  \widetilde{p}_{\varepsilon\delta}\right)  }\right)
T_{\delta}\left(  T_{\varepsilon}\left(  \nabla_{x}\Psi\left(  x\right)
\psi\left(  \frac{x}{\varepsilon}\right)  \phi\left(  \frac{x}{\varepsilon
\delta}\right)  \right.  \right. \\
&  +\left.  \left.  \Psi\left(  x\right)  \nabla_{x}\psi\left(  \frac
{x}{\varepsilon}\right)  \phi\left(  \frac{x}{\varepsilon\delta}\right)
+\Psi\left(  x\right)  \psi\left(  \frac{x}{\varepsilon}\right)
\operatorname{div}_{x}\phi\left(  \frac{x}{\varepsilon\delta}\right)  \right)
\right)  dxdydz.
\end{align*}
Passing to the limit as $\varepsilon$ tends to zero and then using (\ref{5.2})
we obtain%
\begin{align*}
{\int_{\Omega\times Y^{\ast}\times Z^{\ast}}}\widehat{p}\nabla_{x}\Psi\left(
x\right)  \psi\left(  y\right)  \phi\left(  z\right)  dxdydz  &
={\int_{\Omega\times Y^{\ast}\times Z^{\ast}}}\widehat{p}\nabla_{x}\Psi\left(
x\right)  \psi\left(  y\right)  \phi\left(  z\right)  dxdydz-\\
-\int_{\Omega}\widehat{p}\left(  \operatorname{div}_{x}\int_{Y^{\ast}\times
Z^{\ast}}\widehat{u}dydz\right)  dx  &  =-\left\langle \nabla_{x}\widehat
{p},\int_{Y^{\ast}\times Z^{\ast}}\left(  \Psi\left(  x\right)  \psi\left(
y\right)  \phi\left(  z\right)  -\widehat{u}\right)  dydz\right\rangle
_{\Omega}.
\end{align*}
Combining now all the previous convergences we finally get
\begin{align}
&  2\mu{\int_{\Omega\times Y^{\ast}\times Z^{\ast}}}\nabla_{z}\widehat{u}%
\cdot\nabla_{z}\left(  \varPsi-\widehat{u}\right)  dxdydz+g{\int_{\Omega\times
Y^{\ast}\times Z^{\ast}}}\left\vert \nabla_{z}\left(  \varPsi \right)
\right\vert dxdydz-g{\int_{\Omega\times Y^{\ast}\times Z^{\ast}}}\left\vert
\nabla_{z}\widehat{u}\right\vert dxdydz\nonumber\\
&  \geq\left\langle f-\nabla_{x}\widehat{p},\int_{Y^{\ast}\times Z^{\ast}%
}\left(  \varPsi -\widehat{u}\right)  dydz\right\rangle _{\Omega},\ \nonumber
\end{align}
for every $\varPsi\in L^{2}(\Omega,\mathcal{V})$ and by localizing we obtain
(\ref{limprob}).

We notice that the function $\widehat{u}$ which verifies (\ref{limprob}) is
the unique solution of the problem%
\begin{align}
&  2\mu{\int_{ Y^{\ast}\times Z^{\ast}}}\nabla_{z}\widehat{u}\cdot\nabla
_{z}\left(  \widehat{v}-\widehat{u}\right)  dxdydz+g{\int_{ Y^{\ast}\times
Z^{\ast}}}\left\vert \nabla_{z}\widehat{v}\right\vert dxdydz-g{\int_{Y^{\ast
}\times Z^{\ast}}}\left\vert \nabla_{z}\widehat{u}\right\vert
dxdydz\nonumber\\
&  \geq\int_{ Y^{\ast}\times Z^{\ast}}f\left(  \widehat{v}-\widehat{u}\right)
dxdydz,\nonumber
\end{align}
for all $\widehat{v}\in\mathcal{V}$.

The non unique function $\widehat{p}$ corresponding to the pressure is then
recovered by adapting to our case the ideas in \cite{7}. $\blacksquare$

\section{Interpretation of the limit problem\label{section5}}

The limit problem (\ref{limprob}) from Theorem \ref{theor1} can be interpreted
as a non linear Darcy law. In order to derive this result we follow the ideas
in Lions and Sanchez-Palencia \cite{7} for the study of the Bingham flow in a
classical porous medium.

Let $\lambda\in{\mathbb{R}}^{n}$, $v\in\mathcal{V}$ and define
\[
(\lambda,v)_{Y^{\ast}\times Z^{\ast}}=\int_{Y^{\ast}\times Z^{\ast}}%
\lambda_{i}v_{i}dydz.
\]

Denote $\chi(\lambda)=\chi(y,z;\lambda)$ the unique solution of the following
variational inequality: \vskip 2mm Find $\chi(\lambda)\in\mathcal{V}$ such
that%
\begin{align}
&  2\mu%
{\displaystyle\int_{Y^{\ast}}}
{\displaystyle\int_{Z^{\ast}}}
\nabla_{z}\chi(\lambda)\cdot\nabla_{z}\left(  \Phi-\chi(\lambda)\right)
dydz+g%
{\displaystyle\int_{Y^{\ast}}}
{\displaystyle\int_{Z^{\ast}}}
\left\vert \nabla_{z}\left(  \Phi\right)  \right\vert dydz-g%
{\displaystyle\int_{Y^{\ast}}}
{\displaystyle\int_{Z^{\ast}}}
\left\vert \nabla_{z}\chi(\lambda)\right\vert dydz\label{5.1}\\
&  \geq\langle\lambda,\Phi-\chi(\lambda)\rangle_{Y^{\ast}\times Z^{\ast}%
}\nonumber
\end{align}
for every $\Phi\in\mathcal{V}$.

Then we deduce from (\ref{limprob}) and (\ref{5.1}) that%
\[
\widehat{u}(x,y,z)=\chi(y,z;f(x)-\nabla_{x}\widehat{p}(x)).
\]

Relations (\ref{5.2}) and (\ref{5.5}) imply%
\[
\left(  \int_{{Y^{\ast}}\times Z^{\ast}}\widehat{u}(x,y,z)dydz,\nabla
q\right)  _{\Omega} =0,\,\,\,\forall q\in H^{1}(\Omega).
\]
and so the pressure $\widehat{p}$ verifies%
\begin{equation}
\left(  \int_{{Y^{\ast}}\times Z^{\ast}}\chi(y,z;f-\nabla\widehat
{p})dydz,\nabla q\right)  _{\Omega} =0,\,\,\,\forall q\in H^{1}(\Omega).
\label{5.10}%
\end{equation}

Let us now define
\[
\mathcal{K}(\lambda)=\frac{1}{|Y^{\ast}||Z^{\ast}|}\int_{{Y^{\ast}}\times
Z^{\ast}}\chi(y,z;\lambda)dydz,
\]
which is a function from ${\mathbb{R}}^{n}$ into ${\mathbb{R}}^{n}$. Then
relation (\ref{5.10}) reads%
\[
(\mathcal{K}(f-\nabla\widehat{p}),\nabla q)_{\Omega}=0,\,\,\,\forall q\in
H^{1}(\Omega).
\]

Defining the velocity of filtration by
\[
u^{0}(x)=\frac{1}{|Y^{\ast}||Z^{\ast}|}\int_{Y^{\ast}\times Z^{\ast}}%
\widehat{u}(x,y,z)dydz,
\]
we obtain the non linear Darcy law
\begin{equation}
u^{0}(x)=\mathcal{K}(f-\nabla\widehat{p})\,\hbox{ in }\Omega, \label{darcy}%
\end{equation}
where in the right-hand side we have the non linear vectorial function
$\mathcal{K}$.

Moreover, according to (\ref{5.2}) and (\ref{5.5}), function $u^{0}$ verifies%

\[
\hbox{div} u^{0}= 0 \, \hbox { in } \, \Omega,
\]

\[
\nu u^{0} = 0 \, \hbox { on } \, \partial\Omega.
\]

Let us notice that according to Theorem \ref{theor1} we have%

\[
u_{\varepsilon\delta}\rightarrow\frac{\vert Y^{\ast} \vert\vert Z^{\ast}
\vert}{\vert Y \vert\vert Z \vert}u^{0} \hbox { weakly in } (L^{2}\left(
\Omega\right)  )^{n},
\]

and%

\[
\widetilde{p}_{\varepsilon\delta}\rightarrow\widehat{p} \hbox { strongly in }
L_{0}^{2}(\Omega).
\]

This clearly shows that (\ref{darcy}) is the problem verified by the limits of
the sequences $u_{\varepsilon\delta}$ and $\widetilde{p}_{\varepsilon\delta}$,
solutions of (\ref{ineq}).

For seek of completness, we recall below the result obtained for the
homogenization of the Stokes flow in our porous medium and whose limit is a
linear Darcy law. This problem was first studied by Lions in \cite{Lions} with
the method of asymptotic expansions. The justification of the convergence
result is done by Bunoiu and Saint Jean Paulin in \cite{BuSJP}, where the
three-scale convergence method introduced by G. Allaire and M. Briane in
\cite{AllBri} is used.

The Stokes flow can be seen as a particular case of the Bingham flow and it
corresponds to the value zero for the parameter $g$ in the constitutive law.
Indeed, when $g$ equals zero, relation (\ref{sig}) becomes%
\[
\sigma_{ij}=-p_{\varepsilon\delta}\delta_{ij}+2\mu\varepsilon^{2}\delta
^{2}D_{ij}(u_{\varepsilon\delta}),
\]

This particular case corresponds to a Newtonian fluid, which satisfies the
Stokes system:
\begin{align*}
-2\mu\varepsilon^{2}\delta^{2}\Delta u_{\varepsilon\delta} + \nabla
p_{\varepsilon\delta}  &  =f \text{ in }\Omega_{\varepsilon\delta}\\
u_{\varepsilon\delta}  &  =0\text{ on }\partial\Omega_{\varepsilon\delta}.
\end{align*}

In this case, the unique solution $(u_{\varepsilon\delta},p_{\varepsilon
\delta})\in V(\Omega_{\varepsilon\delta})\times L_{0}^{2}(\Omega
_{\varepsilon\delta})$ of the Stokes problem satisfies%
\[
2\mu\varepsilon^{2}\delta^{2}\int_{\Omega_{\varepsilon\delta}}\nabla
u_{\varepsilon\delta}\cdot\nabla vdx=<f-\nabla p_{\varepsilon\delta
},v>_{\Omega_{\varepsilon\delta}},\,\,\forall v\in(H_{0}^{1}(\Omega
_{\varepsilon\delta})^{n}.
\]

Convergence results from Theorem \ref{theor1} for $u_{\varepsilon\delta}$ and
$p_{\varepsilon\delta},$ as far as relations (\ref{5.2})-(\ref{5.9}), hold
true. The only difference is the limit problem (\ref{limprob}) which in this
case reads in the simpler way%
\[
\int_{ Y^{\ast}\times Z^{\ast}}\nabla_{z}\widehat{u}\nabla_{z}\Phi dydz=<
f-\nabla_{x}\widehat{p},\Phi> _{Y^{\ast}\times Z^{\ast}},
\]

for every $\Phi\in\mathcal{V}.$

The linearity of this problem now implies%
\[
\widehat{u}(x,y,z)=\chi(y,z)\left(  f(x)-\nabla_{x}\widehat{p}(x)\right)  ,
\]
where the entries $\chi_{ij},$ $i,j=1,...,n$, of the matrix $\chi$ are the
solutions of the following local problems defined in the domain $Y^{\ast
}\times Z^{\ast}$:

Find $\chi_{i}\in\mathcal{V}$ such that%
\[
\int_{Y^{\ast}\times Z^{\ast}} \nabla_{z} \chi_{i} \nabla_{z} w dydz=(e_{i}%
,w)_{Y^{\ast}\times Z^{\ast}},
\]
for every $w\in\mathcal{V},$ where $e_{i}$ is the $i$-th unit vector of the
canonical base in ${\mathbb{R}}^{n}.$

In this case, the permeability tensor $K$ is defined as the the matrix those
entries are
\[
K_{ij}=\frac{1}{|Y^{\ast}||Z^{\ast}|}\int_{Y^{\ast}\times Z^{\ast}}\chi
_{ij}(y,z)dydz,\ \ \ i,j=1,...,n,
\]
which is linked for every fixed $i=1,...,n$ to the components of the velocity
of filtration via the equality%
\begin{equation}
u_{i}^{0}=K_{ij}\left(  f_{i}-\frac{\partial\widehat{p}}{\partial x_{j}%
}\right)  \,\,\hbox{ in } \Omega, \label{5.11}%
\end{equation}
where we sum over $j$ between $1$ and $n$.

This is the linear Darcy law for our porous medium, which can be also written as%

\[
u^{0}(x)=K(f-\nabla\widehat{p}) \, \hbox { in  } \, \Omega,
\]
where in the right-hand side we multiply a $n\times n $ matrix with a vector
belonging to ${\mathbb{R}}^{n}$.

We observe that the linear Darcy law can be seeen as a particular case of the
non linear one. Indeed, it is obtained when the function $\mathcal{K}$ of
$\lambda$ is linear and so $\mathcal{K}(\lambda)=K\lambda$, where $K$ is a
$n\times n$ matrix.

\end{document}